\documentclass[11pt,reqno]{amsart}
\usepackage{amssymb,amscd,amsbsy}
\usepackage{amssymb,amscd,amsbsy,mathrsfs}
\newcommand{\lb}{\linebreak}

\renewcommand{\a}{\alpha}

\newcommand{\e}{\varepsilon}

\renewcommand{\l}{\lambda}

\newcommand{\s}{\sigma}

\newcommand{\f}{\varphi}

\newcommand{\D}{\Delta}
\renewcommand{\L}{\Lambda}

\renewcommand{\O}{\Omega}

\newcommand{\A}{{\mathscr A}}
\newcommand{\B}{{\mathscr B}}

\newcommand{\E}{{\mathscr E}}

\newcommand{\h}{{\mathscr H}}

\newcommand{\p}{{\mathscr P}}

\newcommand{\X}{{\mathscr X}}
\newcommand{\Y}{{\mathscr Y}}

\newcommand{\C}{{\Bbb C}}

\newcommand{\Z}{{\Bbb Z}}

\newcommand{\bS}{{\boldsymbol S}}

\newcommand{\rf}[1]{(\ref{#1})}

\newcommand{\df}{\stackrel{\mathrm{def}}{=}}

\newcommand{\eeq}{\end{equation}}
\newcommand{\beq}{\begin{equation}}
\newcommand{\bay}{\begin{eqnarray}}
\newcommand{\ba}{\begin{align*}}
\newcommand{\ea}{\end{align*}}
\newcommand{\ey}{\end{eqnarray}}
\newcommand{\bey}{\begin{eqnarray*}}
\newcommand{\eey}{\end{eqnarray*}}

\newcommand{\be}{\infty}

\newcommand{\bl}{\blacksquare}

\newcommand{\Pf}{{\bf Proof. }}

\newcommand{\ov}{\overline}

\newtheorem{thm}{\hspace{\parindent}Theorem}[section]

\newtheorem{cor}[thm]{\hspace{\parindent}Corollary}
\newtheorem{lem}[thm]{\hspace{\parindent}Lemma}

\newtheorem*{rem*}{Remark}

\newcommand\fM{\frak M}

\newcommand{\Bbbone}{{\rm{1\mathchoice{\kern-0.25em}{\kern-0.25em}{\kern-0.2em}{\kern-0.2em}I}}}

\newcommand\mB{\mathcal{B}}

\newcommand\mI{\mathcal{I}}

\begin{document}

\numberwithin{equation}{section}

\numberwithin{equation}{section}

\title{Haagerup tensor products and Schur multipliers}
\author{A.B. Aleksandrov and V.V. Peller}
\thanks{The research on \S\:4 is supported by 
Russian Science Foundation [grant number 23-11-00153].
The research on \S\:5,6,7 is supported by a grant of the Government of the Russian Federation for the state support of scientific research, carried out under the supervision of leading scientists, agreement  075-15-2024-631.
}

%\date{\today}

\

\begin{abstract}
In this paper we compare various classes of Schur multipliers: classical matrix Schur multipliers, discrete Schur multipliers, Schur multipliers with respect to measures and Schur multipliers with respect to spectral measures.
The main result says that in the case of Schur multipliers with respect to measures and spectral measures
such Schur multipliers coincide isometrically with the Haagerup tensor products of the corresponding $L^\be$ spaces. We deduce this result from a well known analogue of it for discrete Schur multipliers.
\end{abstract} 

\maketitle

\section{\bf Introduction}
\setcounter{equation}{0}
\label{In}

\

In this paper we consider various classes of Schur multipliers such as classical matrix Schur multipliers, slightly more general Schur multipliers, Schur multipliers with respect to measures and Schur multipliers with respect to spectral measures.

The notion of a matrix Schur multiplier goes back to the classical paper \cite{Sc}.
We start with the classical case of matrix Schur multipliers. We say that a matrix  
$\{a_{jk}\}_{j,k\ge0}$ is a {\it matrix Schur multiplier} if the following implication holds:
$$
\{b_{jk}\}_{j,k\ge0}\in\mB\quad\Longrightarrow\{a_{jk}b_{jk}\}_{j,k\ge0}\in\mB,
$$
where the symbol $\mB$ denotes the class of infinite matrices that induce bounded linear operators on the space $\ell^2$ of infinite sequences. By the norm of a matrix in $\mB$ we understand 
the operator norm of the corresponding operator on $\ell^2$.

We denote by $\fM$ the pace of matrix Schur multipliers. By the norm  
$\|\{a_{jk}\}_{j,k\ge0}\|_\fM$ we understand the norm of the transformer
\bay
\label{trans}
\{b_{jk}\}_{j,k\ge0}\quad\mapsto\quad\{a_{jk}b_{jk}\}_{j,k\ge0}
\ey
in the operator norm.

It is easy to see from a duality argument that a matrix $\{a_{jk}\}_{j,k\ge0}$ is a Schur multiplier if and only if
it is a Schur multiplier of the trace class, i.e.,
$$
\{b_{jk}\}_{j,k\ge0}\in\bS_1\quad\Longrightarrow\quad\{a_{jk}b_{jk}\}_{j,k\ge0}\in\bS_1,
$$
where the symbol $\bS_1$ denotes the {\it trace class}. Moreover, the norm of the tranformer \rf{trans} in the trace norm coincides with the norm  $\|\{a_{jk}\}_{j,k\ge0}\|_\fM$.

Note also that the norm of a matrix in the space of Schur multipliers can also be defined for finite matrices.

Apparently, the first sufficient condition for an infinite matrix to be a Schur multiplier was obtained in Schur's paper \cite{Sc}: if $\{a_{jk}\}_{j,k\ge0}$ is a matrix of a bounded linear operator on $\ell^2$, then 
$\{a_{jk}\}_{j,k\ge0}$ is a Schur multiplier.

We would like to mention here the paper by G. Bennett \cite{Be}, in which a lot of information on Schur multiplier is gathered.

There are various descriptions of the Schur multipliers. They will be considered below. In particular, 
In the case of matrix Schur multipliers as well as in the more general case of discrete Schur multipliers
there is a known isometric description in terms of the Haagerup tensor product of the spaces 
$\ell^\be$ (see \S\:\ref{discSchur}).

As for Schur multipliers with respect to measures and with repect to spectral measures, there was a known non-isometric description of such Schur multipliers in terms of the Haagerup tensor products of $L^\be$ spaces (see \S\:\ref{Dvoi}).

The principal purpose of this paper is to obtain an {\it isometric description} of the Schur multipliers with respect to measures and with respect to spectral measures in terms the Haagerup tensor products of $L^\be$ spaces.

\

\section{\bf Schur multipliers and  double operator integrals}
\setcounter{equation}{0}
\label{Dvoi}

\

Let us proceed now to a discussion of Schur multipliers with respect to spectral measures. Such multipliers
play an important role in the theory of perturbations of linear operators, see \cite{AP}. To define such Schur multipliers, we need double operator integrals.

Double operator integrals appeared first in the paper \cite{DK} by Yu.L. Daletkii and S.G. Krein.
Later in their papers \cite{BS1,BS2,BS3,BS4}  Birman and Solomyak developed a rigorous mathematical theory of such double operator integrals.

Let $E_1$ and $E_2$ be spectral measures in a separable Hilbert space $\h$. {\it Double operator integrals} are expressions of the form
\bay
\label{dvopin}
\iint\Phi(x,y)\,dE_1(x)T\,dE_2(y).
\ey
Here $\Phi$ is a bounded measurable function and $T$ is a bounded linear operator on Hilbert space.
The approach by Birman and Solomyak  is to define first such double operator integrals in the case when $T$ is an operator of Hilbert--Schmidt class $\bS_2$.
 In this case we can consider the auxiliary spectral measure $\E$ on the Hilbert space $\bS_2$, 
 which can be defined on the measurable rectangles by
$$
\E(\L\times\D)T\df E_1(\L)TE_2(\D),\quad T\in\bS_2.
$$
In \cite{BS4} it was established that the set function defined by this equality extend to a spectral measure. 
In this case the double operator integral \rf{dvopin} is defined in the following way:
$$
\iint\Phi(x,y)\,dE_1(x)T\,dE_2(y)\df\left(\int\Phi\,d\E\right)T.
$$
Clearly, the following inequality holds:
$$
\left\|\iint\Phi(x,y)\,dE_1(x)T\,dE_2(y)\right\|_{\bS_2}\le\|\Phi\|_{L^\be}\|T\|_{\bS_2}.
$$

However, in the case of an arbitrary bounded linear operator $T$,
the double operator integral of the form \rf{dvopin} cannot be defined for arbitrary bounded measurable function $\Phi$, and so it is necessary to impose additional restrictions  on the function.

\medskip

{\bf Definition.} A measurable function $\Phi$ is called a {\it Schur multiplier with respect to spectral measures}
$E_1$ and $E_2$ if the following implication holds: 
$$
T\in\bS_1\quad\Longrightarrow\quad\iint\Phi\,dE_1T\,dE_2\in\bS_1.
$$
We denote the class of all such Schur multipliers by $\fM(E_1,E_2)$. The norm of
 $\Phi$ in $\fM(E_1,E_2)$ is defined as the norm of the operator
 \bay
 \label{Trnaogr}
T\mapsto\iint\Phi\,dE_1T\,dE_2
\ey
on $\bS_1$. If $\Phi\in\fM(E_1,E_2)$, then we can define by duality the transformer 
\rf{Trnaogr} defined on the class of all bounded operators. The norm of this transformer from the class of bounded linear operators on $\h$ into itself coincides by definition with 
$\|\Phi\|_{\fM(E_1,E_2)}$.

\medskip

Consider the case when $\h=\ell^2$ and the spectral measures $E_1$ and $E_2$ on $\h$ are defined by the equality
$$
E_1(\D)u=E_2(\D)u=\sum_{j\in\D}(u,e_j)e_j,\quad\A\in\Z_+,
$$
where $\{e_j\}_{j\ge0}$ is the standard orthonormal basis in $\ell^2$. Then it is easy to observe that 
the function $\Phi$ on $\Z_+\times\Z_+$ belongs to $\fM(E_1,E_2)$ if and only if
the matrix $\{\Phi(j,k)\}_{j,k\ge0}$ is a matrix Schur multiplier.

It i easy to see that if a function $\Phi$ on $\X\times\Y$ belongs to the {\it projective tensor product}
$L^\be(E_1)\hat\otimes L^\be(E_2)$ of the paces $L^\be(E_1)$ and $L^\be(E_2)$ (i.e. $\Phi$ admits a representation
$$
\Phi(x,y)=\sum_{n\ge0}\f_n(x)\psi_n(y),
$$
where $\f_n\in L^\be(E_1)$, $\psi_n\in L^\be(E_2)$ and
$$
\sum_{n\ge0}\|\f_n\|_{L^\be}\|\psi_n\|_{L^\be}<\be),
$$
then $\Phi\in\fM(E_1,E_2)$ and the following equality holds:
$$
\int\limits_\X\int\limits_\Y\Phi(x,y)\,dE_1(x)T\,dE_2(y)=
\sum_{n\ge0}\left(\,\int\limits_\X\f_n\,dE_1\right)T\left(\,\int\limits_\Y\psi_n\,dE_2\right).
$$

More generally, $\Phi\in\fM(E_1,E_2)$ if $\Phi$
belongs to the {\it integral projective tensor product} $L^\be(E_1)\hat\otimes_{\rm i}
L^\be(E_2)$ of the spaces $L^\be(E_1)$ and $L^\be(E_2)$, i.e. $\Phi$ admits a representation
\bay
\label{ipt}
\Phi(x,y)=\int_\O \f(x,w)\psi(y,w)\,d\l(w),
\ey
where $(\O,\l)$ is a space with a $\s$-finite measure, $\f$ i a measurable function on $\X\times \O$,
$\psi$ is a measurable function on $\Y\times \O$ and
\bay
\label{ir}
\int_\O\|\f(\cdot,w)\|_{L^\be(E_1)}\|\psi(\cdot,w)\|_{L^\be(E_2)}\,d\l(w)<\be.
\ey
In this case
$$
\int\limits_\X\int\limits_\Y\Phi(x,y)\,dE_1(x)T\,dE_2(y)=
\int\limits_\O\left(\,\int\limits_\X\f(x,w)\,dE_1(x)\right)T
\left(\,\int\limits_\Y\psi(y,w)\,dE_2(y)\right)\,d\l(w).
$$
Clearly, the function
$$
w\mapsto \left(\,\int_\X\f(x,w)\,dE_1(x)\right)T
\left(\,\int_\Y\psi(y,w)\,dE_2(y)\right)
$$
is weakly measurable and
$$
\int\limits_\O\left\|\left(\,\int\limits_\X\f(x,w)\,dE_1(x)\right)T
\left(\,\int\limits_\Y\psi(y,w)\,dE_2(w)\right)\right\|\,d\l(w)<\be.
$$
Moreover,
$$
\|\Phi\|_{\fM(E_1,E_2)}\le\|\Phi\|_{L^\be(E_1)\hat\otimes_{\rm i}L^\be(E_2)},
$$
where $\|\Phi\|_{L^\be\hat\otimes_{\rm i}L^\be}$ is the infimum of the integral on the left-hand side of 
\rf{ir} over all representations of $\Phi$ in the form \rf{ipt}. 

It turns out that all Schur multipliers can be obtained this way, see \cite{Pe}.

Another sufficient condition for the membership in $\fM(E_1,E_2)$ can be stated in terms of the  {\it Haagerup tensor product}  
$L^\be(E_1)\!\otimes_{\rm h}\!L^\be(E_2)$, whech can be defined a the space of functions $\Phi$ of the form
\bay
\label{FiH}
\Phi(x,y)=\sum_{n\ge0}\f_n(x)\psi_n(y),
\ey
where $\f_n\in L^\be(E_1)$, $\psi_n\in L^\be(E_2)$; moreover,
$$
\{\f_n\}_{n\ge0}\in L_{E_1}^\be(\ell^2)\quad\mbox{and}\quad
\{\psi_n\}_{n\ge0}\in L_{E_2}^\be(\ell^2).
$$
By the {\it norm of a function $\Phi$ in $L^\be(E_1)\!\otimes_{\rm h}\!L^\be(E_2)$} we understand the infimum of the expression
\bay
\label{normaHaag}
\big\|\{\f_n\}_{n\ge0}\big\|_{L_{E_1}^\be(\ell^2)}
\big\|\{\psi_n\}_{n\ge0}\big\|_{L_{E_2}^\be(\ell^2)}
\ey
over all representations of $\Phi$ of the form \rf{FiH}. Here
$$
\big\|\{\f_n\}_{n\ge0}\big\|_{L_{E_1}^\be(\ell^2)}\df
\Big\|\sum_{n\ge0}|\f_n|^2\Big\|_{L^\be(E_1)}^{1/2}\quad\!\!\mbox{and}\quad\!\!
\big\|\{\psi_n\}_{n\ge0}\big\|_{L_{E_1}^\be(\ell^2)}\df
\Big\|\sum_{n\ge0}|\psi_n|^2\Big\|_{L^\be(E_2)}^{1/2}.
$$
It is easy to verify that if $\Phi\in L^\be(E_1)\!\otimes_{\rm h}\!L^\be(E_2)$, then $\Phi\in\fM(E_1,E_2)$ and
$$
\iint\Phi(x,y)\,dE_1(x)T\,dE_2(y)=
\sum_{n\ge0}\Big(\int\f_n\,dE_1\Big)T\Big(\int\psi_n\,dE_2\Big).
$$
The series on the right converges in the weak operator topology and
$$
\|\Phi\|_{\fM(E_1,E_2)}\le\|\Phi\|_{L^\be(E_1)\otimes_{\rm h}L^\be(E_2)}.
$$

As one can see from the following assertion, the condition $\Phi\in L^\be(E_1)\!\otimes_{\rm h}\!L^\be(E_2)$ 
is not only sufficient but also necessary. 

\medskip

{\bf A criterion of the membership in the space of Schur multipliers.}
{\it Let $\Phi$ be a measurable function on
$\X\times\Y$ and let $\mu$ and $\nu$ be positive measure on $\X$ and $\Y$ that are mutually absolutely continuous with respect to  $E_1$ and $E_2$. The following statements are equivalent:

{\rm (i)} $\Phi\in\fM(E_1,E_2)$;

{\rm (ii)} $\Phi\in L^\be(E_1)\hat\otimes_{\rm i}L^\be(E_2)$;

{\rm (iii)} $\Phi\in L^\be(E_1)\!\otimes_{\rm h}\!L^\be(E_2)$;

{\rm (iv)} there exist a $\s$-finite measure space $(\O,\l)$, 
measurable functions $\f$ on $\X\times\O$ and $\psi$ on $\Y\times\O$ such that 
{\em\rf{ipt}} holds and
\bay
\label{bs}
\left\|\left(\int_\O|\f(\cdot,w)|^2\,d\l(w)\right)^{1/2}\right\|_{L^\be(E_1)}
\left\|\left(\int_\O|\psi(\cdot,w)|^2\,d\l(w)\right)^{1/2}\right\|_{L^\be(E_2)}<\be;
\ey

{\rm (v)} If the integral operator $\mI_h$
from $L^2(\nu)$ to $L^2(\mu)$ belongs to $\bS_1$, then the same can be said
about the integral operator
$\mI_{\Phi h}$,
where
$$
(\mI_h f)(y)\df\int h(x,y)f(y)\,d\nu(y).
$$
Moreover, $\|\Phi\|_{\fM(E_1,E_2)}$ coincides with the norm of the transformer 
$$
\mI_h\mapsto \mI_{\Phi h}
$$
in the trace norm.
}

\medskip

We refer the reader to the paper \cite{Pe}, in which a proof of the above assertion is essentially given
as well as to \cite{AP2} for the proof of the equivalence of (i) and (v).

It is well known (see \cite{Pi}) that in the case of discrete Schur multipliers the space of Schur multipliers
coincides with the Haagerup tensor product $\ell^\be\!\otimes_{\rm h}\ell^\be$ isometrically, see \S\:\ref{discSchur}).

The main result of this paper is that in the case of Schur multipliers with respect to spectral meaures
$E_1$ and $E_2$, the space of Schur multipliers $\fM(E_1,E_2)$ coincides with the Haagerup tensor product
$L^\be(E_1)\!\otimes_{\rm h}\!L^\be(E_2)$ also isometrically.

\

\section{\bf The role of discrete Schur multipliers}
\setcounter{equation}{0}
\label{discSchur}

\

Let $\h_1$ and $\h_2$ be Hilbert spaces. we denote by $\mathcal B(\h_1,\h_2)$ the space of bounded linear operators from $\h_1$ to $\h_2$. The space of trace class operators from $\h_1$ to $\h_2$ is denoted by
$\bS_1(\h_1,\h_2)$.

Let $S$ and $T$  be arbitrary set. A function $\Phi : S\times T\to\C$ is called a {\it Schur multiplier on the space} $\mathcal B(\ell^2(T),\ell^2(S))$ if for an arbitrary operator 
$A\in \mathcal B(\ell^2(T),\ell^2(S))$ with matrix $\{a(s,t)\}_{(s,t)\in S\times T}$, the matrix 
$\Phi*a\df\{\Phi(s,t)a(s,t)\}_{(s,t)\in S\times T}$ also induces a bounded linear operator
from $\ell^2(T)$ to $\ell^2(S)$, which will be denoted by $M_\Phi A$. 

Let $\frak M(S\times T)$ denote the space of Schur multipliers on 
$\mathcal B(\ell^2(T),\ell^2(S))$. Let $\Phi\in\frak M(S\times T)$. 
It is easy to see that the transformer $M_\Phi$, which associates with the operator $A\in \mathcal B(\ell^2(T),\ell^2(S))$ 
the operator $M_\Phi A\in \mathcal B(\ell^2(T),\ell^2(S))$ is a continuous operator on the space $\mathcal B(\ell^2(T),\ell^2(S))$.
We denote by $\|\Phi\|_{\frak M(S\times T)}$ the norm of the transformer $M_\Phi$.

It is clear from duality arguments that we obtain an equivalent definition of $\frak M(S\times T)$ if instead of the space $\mathcal B(\ell^2(T),\ell^2(S))$ we consider the space
$\bS_1(\ell^2(T),\ell^2(S))$. In other words, a function $\Phi : S\times T\to\C$ is a Schur multiplier on  
$\mathcal B(\ell^2(T),\ell^2(S))$ if and only if for an arbitrary operator $A$ of class $\bS_1(\ell^2(T),\ell^2(S))$, the operator $M_\Phi A$ is also of trace class; moreover,
\bey
\|\Phi\|_{\frak M(S\times T)}=\sup\{\|M_\Phi A\|_{\bS_1(\ell^2(T),\ell^2(S))}: A\in\bS_1(\ell^2(T),\ell^2(S)), \|A\|_{\bS_1(\ell^2(T),\ell^2(S))}\le1\}\\
=\sup\{\|M_\Phi A\|_{\mathcal B(\ell^2(T),\ell^2(S))}: A\in\mathcal B(\ell^2(T),\ell^2(S)), \|A\|_{\mathcal B(\ell^2(T),\ell^2(S))}\le1\}.
\eey

\begin{thm}  
\label{pmt0}
Let $\Phi : S\times T\to\C$ be a function and let $C$ be a positive number.
The following are equivalent: 

{\em(i)} $\Phi\in\frak M(S\times T)$ and $\|\Phi\|_{\frak M(S\times T)}\le C$,

{\em(ii)} there exist families of vectors $\{x_s\}_{s\in S}$ and $\{y_t\}_{t\in T}$ in a (not necessarily separable) Hilbert space such that such that $\Phi(s,t)=\langle x_s,y_t\rangle$ for all $(s,t)\in S\times T$ and 
$$
\sup\limits_{s\in S}\|x_s\|\sup\limits_{t\in T}\|y_t\|\le C.
$$
\end{thm}

This theorem is given in Theorem 5.1 of Pisier's monograph \cite{Pi}, which contains a number of historical remarks on this result.

Note that the implication (ii)$\Longrightarrow$(i) is absolutely elementary.

We are going to make a simple remark to this theorem, a similar remark 
 Remark is made to Theorem 2.2.2 in the survey \cite{AP}.

\medskip

{\bf Remark.} 
\label{plot0}
We can also require in  Theorem \ref{pmt0} that the linear span of each of the families  $\{x_s\}_{s\in S}$ and
$\{ y_t\}_{t\in T}$ be dense in that Hilbert space.

\medskip

In this paper we are going to use abolutely elementary methods to deduce from Theorem \ref{pmt0} an analog of it for other kinds of Schur multipliers.

\

\section{\bf Schur multiplier with respect to measures}

\

Let $\frak A$ be a $\s$-algebra of subsets of $S$ and let $\frak B$ be  a $\s$-algebra of subsets of  $T$.
Suppose that $\mu$ and $\nu$ are $\s$-finite measures on the $\s$-algebras $\frak A$ and $\frak B$.

We are going to show that Theorem  \ref{pmt0} allows one to obtain an analogue of it for Schur multipliers with respect to measures.

With each function $h$ in $L^2(\mu\times\nu)$ we associate the integral operator 
$\mathcal I_h:L^2(\nu)\to L^2(\mu)$,
$$
(\mathcal I_h f)(s)=\int_T h(s,t) f(t)\,d\nu(t),\quad f\in L^2(\nu).
$$

A function $\Phi$ on the set $S\times T$ is called {\it a Schur multiplier with respect to measures $\mu$ and 
$\nu$} if it is $\mu\times\nu$-measurable and the following implication holds:
$$
\mI_h\in\bS_1\quad\Longrightarrow\quad\mI_{\Phi h}\in\bS_1.
$$
In this case by the closed graph theorem, there exists a positive number $c$ such that
\bay
\label{pkb}
\|\mathcal I_{\Phi h}\|_{\bS_1(L^2(\nu),L^2(\mu))}\le c\|\mathcal I_h\|_{\bS_1(L^2(\nu),L^2(\mu))}
\ey
for all $h$ in $L^2(\mu\times\nu)$ (We assume that $\|\mathcal I_h\|_{\bS_1(L^2(T,\nu),L^2(S,\mu))}=\be$,
if the oprrator $\mathcal I_h$ is not of trace class).

Let $\frak M(\mu,\nu)$ denote the set of all such Schur multipliers. 
For $\Phi\in \frak M(\mu,\nu)$, we denote by
denote by $\|\Phi\|_{\frak M(\mu,\nu)}$ the smallest constant $c\ge0$ satisfying
\rf{pkb}. 

It is well known and it is easy to verify that
 $\frak M(\mu,\nu)\subset L^\be(\mu\times\nu)$ and 
$\|\Phi\|_{L^\be(\mu\times\nu)}\le \|\Phi\|_{\frak M(\mu,\nu)}$, $\Phi\in\frak M(\mu,\nu)$.

It is easy to verify by duality that $\Phi\in\frak M(\mu,\nu)$ if and only if
$$
\|\mathcal I_{\Phi h}\|_{\mathcal B(L^2(\nu),L^2(\mu))}\le\|\Phi\|_{\frak M(\mu,\nu)}
\|\mathcal I_h\|_{\mathcal B(L^2(\nu),L^2(\mu))}
$$
for an arbitrary function $h$ in $L^2(\mu\times\nu)$. Moreover,  $\|\Phi\|_{\frak M(\mu,\nu)}$ is the smallest constant $c\ge0$, for which the inequality
$$
\|\mathcal I_{\Phi h}\|_{\mathcal B(L^2(\nu),L^2(\mu))}\le c\|\mathcal I_h\|_{\mathcal B(L^2(\nu),L^2(\mu))}
$$
holds for all $h\in L^2(\mu\times\nu)$.

\medskip

%\begin{lem}  Пусть $g\in L^0(\mu)$ и $h\in L^0(\nu)$. Предположим, что функции $g$ и $h$ 
%положительны $\mu$-почти всюду на $S$ и $\nu$-почти всюду на $T$  соответственно. 
%Рассмотрим меры $\mu_0$
%и $\nu_0$ такие, что $d\mu_0=g d\mu$ и $d\nu_0=h d\nu$.  Пусть $K\in L^2(\mu\times\nu)$.
%Положим $K_0(s,t)=\frac{K(s,t)}{\sqrt{g(s)h(t)}}$. Тогда
%$\|\mathcal I_{K_0}\|_{\mathcal B(L^2(\nu_0),L^2(\mu_0))}=\|\mathcal I_K\|_{\mathcal B(L^2(\nu),L^2(\mu))}$.
%\end{lem}

The following assertion is well known. We give it here for completeness.

\begin{lem}  
Let $\xi$ be a positive $\mu$-measurable function on $S$ and let
$\eta$ be a positive $\nu$-measurable function on $T$.
We define the measures $\mu_0$
and $\nu_0$ by $d\mu_0=\xi d\mu$ and $d\nu_0=\eta d\nu$.  Given $h\in L^2(\mu\times\nu)$, we put
$h_0(s,t)=\frac{h(s,t)}{\sqrt{\xi(s)\eta(t)}}$. Then
$\|\mathcal I_{h_0}\|_{\mathcal B(L^2(\nu_0),L^2(\mu_0))}=\|\mathcal I_h\|_{\mathcal B(L^2(\nu),L^2(\mu))}$.
\end{lem}

\Pf Consider the operators $U_\xi: L^2(\mu)\to L^2(\mu_0)$ and $V_\eta: L^2(\nu)\to L^2(\nu_0)$ defined by 
$U_\xi f=\frac f{\sqrt\f}$ and $V_\eta f=\frac f{\sqrt\psi}$. Clearly, $U_\xi$ and $V_\eta$ are unitary operators.
It remains to observe that  $\mathcal I_{h_0}=U_\xi \mathcal I_h V_\eta^{-1}$. $\bl$

Therefore, under the hypotheses of this lemma, the following equalities hold:
$$
\frak M(\mu_0,\nu_0)=\frak M(\mu,\nu)\quad\text{and}\quad  \|\cdot\|_{\frak M(\mu_0,\nu_0)}=\|\cdot\|_{\frak M(\mu,\nu)}.
$$

\medskip

Therefore when studying the Schur multipliers $\frak M(\mu,\nu)$m the case of $\s$-finite measures
$\mu$ and $\nu$ reduces to the case of finite measures $\mu$ and $\nu$.

 \medskip

%{\bf Замечание 2.}  Пусть $g\in L^0(\mu)$ и $h\in L^0(\mu)$. Предположим, что функции $g$ и $h$ 
%положительны почти всюду на $S$ и на $T$  соответственно. Рассмотрим меры $\mu_0$
%и $\nu_0$ такие, что $d\mu_0=g d\mu$ и $d\nu_0=h d\nu$.  Легко видеть, что
%$\frak M(\mu_0,\nu_0)=\frak M(\mu,\nu)$ и $\|\Phi\|_{\frak M(\mu_0,\nu_0)}=\|\Phi\|_{\frak M(\mu,\nu)}$
%для всех $\Phi\in\frak M(\mu,\nu)$.
%
%
%\medskip

In the statement of the next theorem we use the notation
$L^p(\mu,X)$ for the $L^\p$ space of 
$X$-valued
$\mu$-measurable functions, where $X$ is a Banach space .

\begin{thm}  
\label{pmtn}
Let $(S,\frak A,\mu)$ and $(T,\frak B,\nu)$ be spaces with $\s$-finite measures,
\lb$\Phi\in L^\be(\mu\times\nu)$ and let $C$ be a positive integer.
The following are equivalent:

{\em(i)} $\Phi\in\frak M(\mu,\nu)$  and $\|\Phi\|_{\frak M(\mu,\nu)}\le C$;

{\em(ii)} there exist a Hilbert space $\h$, 
functions $x\in L^\be(\mu,\mathscr H)$ and $y\in L^\be(\nu,\mathscr H)$,
such that
$\Phi(s,t)=\langle x(s),y(t)\rangle$ for $\mu\times\nu$-almost all $(s,t)\in S\times T$ and
$$
\|x\|_{L^\be(\mu,\mathscr H)}\|y\|_{L^\be(\nu,\mathscr H)}\le C;
$$

{\em(iii)} $\Phi$ can be represented in the form
$$
\Phi(s,t)=\sum_{n=1}^\be\f_n(s)\psi_n(t)\quad\mu\times\nu\text{-almost everywhere on}\quad S\times T
$$
so that
$$
\sum_{n=1}^\be|\f_n|^2\le C\quad\mu\text{-almost everywhere and}\quad
\sum_{n=1}^\be|\psi_n|^2\le C\quad\nu\text{-almost everywhere},
$$
where  $\f_n\in L^\be(\mu)$ and $\psi_n\in L^\be(\nu)$ for all $n\ge1$.
\end{thm}

{\bf The proof of simplest implications.}  Let us first prove that (ii)$\Longleftrightarrow$(iii).
Suppose that (ii) holds. To deduce (iii), it suffice to select an orthonormal basis $\{e_n\}_{n=1}^\be$ in 
$\mathscr H$ and put $\f_n(s)\df\langle x(s),e_n\rangle$,
$\psi_n(t)\df\langle e_n,y(t)\rangle$. In the cae of a finite-dimensional space $\mathscr H$ this reasoning
gives us finite sequences $\{\f_n(s)\}_{n=1}^N$ and $\{\psi_n(t)\}_{n=1}^N$.
In this case we can put $\f_n\df 0$ and $\psi_n\df 0$ for $n>N$.

To deduce (ii) from (iii), it suffice to put $\mathscr H\df\ell^2$, 
$x(s)\df\{\f_n(s)\}_{n=1}^\be$ and $y(t)\df\{\ov{\psi_n(t)}\}_{n=1}^\be$.

Let us prove now that, (ii)$\Longrightarrow$(i).

Suppose that (ii) holds. We have to establish the following inequality:
$$
\|\mathcal I_{\Phi h}\|_{\bS_1(L^2(T,\nu),L^2(S,\mu))}\le C\|\mathcal I_h\|_{\bS_1(L^2(T,\nu),L^2(S,\mu))}
$$
for all $h$ in $L^2(\mu\times\nu)$. It is easy to see that it suffice to prove this inequality only in the case when
$I_h$ has rank $1$. In this case $h$ can be represented as $h(s,t)=u(s)\ov{v(t)}$ with 
$u\in L^2(\mu)$ and $v\in L^2(\nu)$.
Then $\Phi(s,t)h(s,t)=\langle u(s)x(s),v(t)y(t)\rangle$. Therefore,
\bey
\|\mathcal I_{\Phi h}\|_{\bS_1(L^2(T,\nu),L^2(S,\mu))}=\|ux\|_{L^2(\mu)}\|vy\|_{L^2(\nu)}
\le\|g\|_{L^2(\mu)}\|x\|_{L^\be(\mu,\mathscr H)}\|h\|_{L^2(\nu)}\|y\|_{L^\be(\nu,\mathscr H)}\\
\le C\|u\|_{L^2(\mu)}\|v\|_{L^2(\nu)}=C\|\mathcal I_h\|_{\bS_1(L^2(T,\nu),L^2(S,\mu))}.\quad\bl
\eey

The following lemma will play an important role in the proof of the implication (i)$\Longrightarrow$(iii) in Theorem \ref{pmtn}.

\begin{lem} 
\label{STmunu}
Let $\Phi\in\frak M(S\times T)$. Suppose that $\Phi$ is measurable with respect to the measure
$\mu\times\nu$ on $S\times T$. Then $\Phi$ can be represented in the form
$$
\Phi(s,t)=\sum_{n=1}^\be\f_n(s)\psi_n(t)\quad\mu\times\nu\mbox{-almost~everywhere~on}\quad S\times T
$$
so that
$$
\sum_{n=1}^\be|\f_n|^2\le\|\Phi\|_{\frak M(S\times T)}\quad\mu\mbox{-almost~everywhere}
$$
and
$$
\sum_{n=1}^\be|\psi_n|^2\le\|\Phi\|_{\frak M(S\times T)}\quad\nu\mbox{-almost~everywhere},
$$
where  $\f_n\in L^\be(\mu)$ and $\psi_n\in L^\be(\nu)$ for all $n\ge1$.
\end{lem}

\Pf Clearly, $\Phi$ is bounded. Thus the fact that $\Phi$ is $\mu\times\nu$-measurable implies that  for 
$\mu$-almost all $s$ in $S$, the function $\Phi(s,\cdot)$ belong to $L^\be(T,\nu)$ 
and for $\nu$-almost all $t$ in $T$, the function $\Phi(\cdot,t)$ belong to $L^\be(S,\mu)$.

Let us first prove the theorem in the case when for all $s$ in $S$, the function
$\Phi(s,\cdot)$ belongs to $L^\be(T,\nu)$ and for all $t\in T$, the function
$\Phi(\cdot,t)$ belongs to $L^\be(S,\mu)$.

It follows from Theorem \ref{pmt0} that $\Phi$ can be represented as $\Phi(s,t)=\langle\f(s),\psi(t)\rangle$
for functions $\f:S\to\mathscr H$ and $\psi:T\to\mathscr H$ so that $\|\f(s)\|_{\mathscr H}^2\le\|\Phi\|_{\frak M(S\times T)}$ for all $s\in S$ and
$\|\psi(t)\|_{\mathscr H}^2\le\|\Phi\|_{\frak M(S\times T)}$ for all $t\in T$.

It follows from the Remark to Theorem
\ref{pmt0} that we can also select $\f$ and $\psi$ so that the linear span of each of the families
 $\{\f(s)\}_{s\in S}$ and $\{\psi(t)\}_{t\in T}$ be dense in $\mathscr H$. %$\f$ и $\psi$ слабо измеримы на $S$ и $T$ соответственно. 
Then $\f$ is weakly $\mu$-measurable while $\psi$ is weakly $\nu$-measurable. 

Let $\{e_\a\}_{\a\in\A}$ be an orthnormal basis in $\mathscr H$.
Put $\f_\a(s)\df\langle\f(s),e_\a\rangle$ and $\psi_\a(t)\df\langle e_\a,\psi(t)\rangle$. Clearly, for all
 $\a$ in $\A$, the functions $\f_\a$ and $\psi_\a$ are $\mu$-measurable and $\nu$-measurable
respectively. The equality $\Phi(s,t)=\langle\f(s),\psi(t)\rangle$ can be rewritten in the following way:
$$
\Phi(s,t)=\sum_{\a\in \A}\f_\a(s)\psi_\a(t)\quad\text{for all}\quad (s,t)\in S\times T;
$$
moreover,
$$
\sum_{\a\in\A}|\f_\a(s)|^2\le\|\Phi\|_{\frak M(S\times T)}\quad\mbox{for all}\quad s\in S\quad
$$
and
$$
\sum_{\a\in\A}|\psi_\a(t)|^2\le\|\Phi\|_{\frak M(S\times T)}
\,\text{for all}\, t\in T.
$$

It follows that
$$
\sum_{\a\in\A}\|\f_\a\|_{L^2(\mu)}^2\le\|\Phi\|_{\frak M(S\times T)}\quad\text{and}\quad
\sum_{\a\in\A}\|\psi_\a\|_{L^2(\nu)}^2\le\|\Phi\|_{\frak M(S\times T)}.
$$
Indeed, this is obvious if we replace the summation over $\A$ with the summation over an arbitrary at most countable subset $\B$ of $\A$. It remains to pass to the supremum over all such subsets $\B$.

Let us show that there exists an at most countable subset  $\A_0$ of $\A$ 
such that
%\bay
%\label{rav0}
%\Phi(s,t)=\sum_{\a\in \A_0}\f_\a(s)\psi_\a(t)\quad\text{почти всюду на}\quad  S\times T.
%\ey
\bay
\label{rav0}
\Phi(s,t)=\sum_{\a\in \A_0}\f_\a(s)\psi_\a(t)\quad\,\,\mu\times\nu\text{-almost everywhere on}\quad  S\times T.
\ey
We can select $\A_0$ to be the set of $\a\in\A$ such that
$\|\f_\a\|_{L^2(S)}>0$. Obviouly, $\A_0$ is at most countable and 
$\f_\a=0$ almost everywhere on $S$ for all $\a$ in $\A\setminus\A_0$. 

To establish \rf{rav0}, it suffices to show that
$$
\int_T\left(\int_S\Big|\Phi(s,t)-\sum_{\a\in\A_0}\f_\a(s)\psi_\a(t)\Big|\,d\mu(s)\right)\,d\mu(t)=0.
$$
%$$
%\int_S\Big|\sum_{\a\in\A\setminus \A_0}\f_\a(s)\ov{\psi_\a(t)}\Big|\,d\mu(s)=0.
%$$
We fix a point $t$ in $T$. Let $\A(t)$ denotes the set of $\a$ in $\A\setminus \A_0$ such that $\psi_\a(t)\ne0$.
Clearly, $\A(t)$ is at most countable and
\bey
\int_S\Big|\Phi(s,t)-\sum_{\a\in\A_0}\f_\a(s)\psi_\a(t)\big|\,d\mu(s)\le\int_S\sum_{\a\in\A(t)}|\f_\a(s)|\cdot|\psi_\a(t)|\,d\mu(s)\\
=\sum_{\a\in\A(t)}|\psi_\a(t)|\int_S|\f_\a(s)|\,d\mu(s)=0
\eey
for all $t$ in $T$. 

%Теперь ясно, что из доказанной импликации (ii)$\Longrightarrow$(i) теоремы \ref{pmtn} следует, что
%$\Phi\in\frak M(\mu,\nu)$ и $\|\Phi\|_{\frak M(\mu,\nu)}\le\|\Phi\|_{\frak M(S\times T)}$.

Let us consider now the general case. Clearly, there exists a $\mu$-measurable subset $S_0$ of $S$ and a 
$\nu$-measurable subset $T_0$ of $T$ such that
$\mu(S\setminus S_0)=0$, $\nu(T\setminus T_0)=0$, the function 
$\Phi(s,\cdot)$ belongs to $L^\be(T,\nu)$ for all $s$ in $S_0$ and the function $\Phi(\cdot,t)$ belongs to $L^\be(S,\mu)$ for all $t$ in $T_0$.
 
Put $\widetilde\Phi\df\chi_{S_0\times T_0}\Phi$. It remains to observe that  for the function
$\widetilde\Phi$ the theorem has already been proved and 
$\|\widetilde\Phi\|_{\frak M(S\times T)}\le \|\Phi\|_{\frak M(S\times T)}$. $\bl$

\begin{cor} 
\label{STmunu0}
Let $\Phi\in\frak M(S\times T)$. Suppose that the function $\Phi$  is measurable with respect to the measure
$\mu\times\nu$ on $S\times T$.
Then $\Phi\in\frak M(\mu,\nu)$ and $\|\Phi\|_{\frak M(\mu,\nu)}\le\|\Phi\|_{\frak M(S\times T)}$
\end{cor}

\Pf This  result follows from  the lemma if we use the implication (iii)$\Longrightarrow$(i) of Theorem \ref{pmtn} that has already been proved. $\bl$

\

\section{\bf Auxiliary results}
\label{pomoshch'}

\

Let $(S,\frak A,\mu)$ and $(T,\frak B,\nu)$ be spaces with finite measures. Suppose that
$\{S_j\}_{j=1}^m$ and $\{T_k\}_{k=1}^n$ are are finite measurable partitions of the sets $S$ and $T$ such that
$\mu(S_j)>0$ for all $j$ and $\nu(T_k)>0$ for all $k$.

With each function $f$ in $L^2(\mu)$ we associate the function $\mathcal P f$ whose value at each point of  $S_j$ is the mean value of $f$  over $S_j$ for all $j$. Clearly,  $\mathcal P$ is an orthogonal projection on   $L^2(\mu)$ and
 $$
 (\mathcal P f)(x)=\int_S\left(\sum_{j=1}^m\frac1{\mu(S_j)}\Bbbone_{S_j}(x)\Bbbone_{S_j}(y)\right)f(y)\, d\mu(y).
 $$
 We also consider a similar orthogonal projection
  $\mathcal Q$ on the space $L^2(\nu)$,
 $$
 (\mathcal Q f)(x)=\int_T\left(\sum_{k=1}^n\frac1{\nu(T_k)}\Bbbone_{T_k}(x)\Bbbone_{T_k}(y)\right)f(y)\, d\nu(y).
 $$
 
%{\bf Вычисления для меня.}
% 
%\bey
%(\mathcal I_K\mathcal Q f)(s)=\int_T K(s,t) (\mathcal Qf)(t)\,d\nu(t)\\
%=\int_T\int_T K(s,t)\sum_{k=1}^n\left(\frac1{\nu(T_k)}\Bbbone_{T_k}(t)\Bbbone_{T_k}(y)\right)f(y)\, d\nu(t) d\nu(y)
%\eey
%
%\bey
%(\mathcal PI_K\mathcal Q f)(x)\\
%=\int_S\int_T\int_T K(s,t)\sum_{k=1}^n\left(\frac1{\nu(T_k)}\Bbbone_{T_k}(t)\Bbbone_{T_k}(y)\right)
%\sum_{j=1}^m\left(\frac1{\mu(S_j)}\Bbbone_{S_j}(x)\Bbbone_{S_j}(s)\right)f(y)\, d\nu(t) d\nu(y) d\nu(s)
%\eey

We can also consider the orthogonal projection $\mathcal P\odot\mathcal Q$ 
on the space $L^2(\mu\times\nu)$ associated with the finite partition  $\{S_j\times T_k\}_{1\le j\le m,1\le k\le n}$
of the set  $S\times T$.
Let $h\in L^2(\mu\times\nu)$.  It is easy to see that  
$\mathcal I_{(\mathcal P\odot\mathcal Q)h}=\mathcal P \mathcal I_h\mathcal Q$.
It follows that
$\|\mathcal I_{(\mathcal P\odot\mathcal Q)h}\|_{\mathcal B(L^2(\nu),L^2(\mu))}
\le\|\mathcal I_h\|_{\mathcal B(L^2(\nu),L^2(\mu))}$.
Thus, we have proved the following lemma.

\begin{lem}  
\label{PQK0}
Let $h\in L^2(\mu\times\nu)$. Then
$\|\mathcal I_{(\mathcal P\odot\mathcal Q)h}\|_{\mathcal B(L^2(\nu),L^2(\mu))}\le\|
\mathcal I_h\|_{\mathcal B(L^2(\nu),L^2(\mu))}$.
\end{lem}

 \begin{lem}  
\label{nmatr0}
Let $C=\{c_{jk}\}$ be a matrix of size $m\times n$.  Put  $h(s,t)=\lb(\mu(S_j))^{-\frac12}c_{jk}(\nu(T_k))^{-\frac12}$
for $(s,t)\in S_j\times T_k$. Then $\|C\|=\|\mathcal I_h\|_{\mathcal B(L^2(\nu),L^2(\mu))}$.
\end{lem}

\Pf It follow from the equality $\mathcal I_h=\mathcal I_h\mathcal Q$ that the norm of the operator $\mathcal I_h$ is attained at a function that i constant on each set $T_k$ ($1\le k\le n$). This allows u to observe that the norm of the operator $\mathcal I_h$ is equal o the norm of the matrix $C$. $\bl$

%\begin{lem}  
%\label{mmatr0}
%Пусть $C=\{c_{jk}\}$ матрица размеров $m\times n$.  Положим
%$\Phi(s,t)=c_{jk}$
%при $(s,t)\in S_j\times T_k$. Тогда $\|C\|_{\frak M}=\|\Phi\|_{\frak M(S\times T)}=\|\Phi\|_{\frak M(\mu,\nu)}$.
%\end{lem}

\begin{lem}  
\label{mmatr0}
Let $C=\{c_{jk}\}$ be a matrix of size $m\times n$.  Put
$\Phi(s,t)=c_{jk}$
for $(s,t)\in S_j\times T_k$. Then $\|C\|_{\frak M}=\|\Phi\|_{\frak M(\mu,\nu)}$.
\end{lem}

\Pf The inequality  $\|C\|_{\frak M}\le\|\Phi\|_{\frak M(\mu,\nu)}$
is obvious. The opposite \lb inequality follows easily from the implication (ii)$\Longrightarrow$(i) in Theorem \ref{pmtn} that has already been proved. $\bl$

Lemmata \ref{nmatr0} and \ref{mmatr0} imply the following assertion.

\begin{lem}  
\label{msup0}
Let $\Phi\in L^2(\mu\times\nu)$.
Suppose that $(\mathcal P\odot\mathcal Q)\Phi=\Phi$ almost everywhere with respect to the measure $\mu\times\nu$.
Then 
$$
\|\Phi\|_{\frak M(\mu,\nu)}=\sup\{\|\mathcal I_{\Phi h}\|_{\mathcal B(L^2(\nu),L^2(\mu))}\},
$$
where the supremum is taken over all $h\in L^2(\mu\times\nu)$ such that $(\mathcal P\odot\mathcal Q)h=h$
and $\|\mathcal I_h\|_{\mathcal B(L^2(\nu),L^2(\mu))}\le1$.
\end{lem}

%\begin{lem} Пусть $(\mS,\frak A,\mu)$ и $(\mT,\frak B,\nu)$ пространства c  конечными мерами.
%Пусть $\{\mS_j\}_{j=1}^m$ и  $\{\mT_k\}_{k=1}^n$ -- конечные измеримые разбиения. Предположим,
%что $\mu(\mS_j)>0$ для всех $j$ и $\nu(\mT_k)>0$ для всех $k$. Пусть $\Phi\in\frak M(\mu,\nu)$.
%Положим $a_{jk}$ обозначает среднее значение функции $\Phi$  на множестве $\mS_j\times \mT_k$,
%$$
%c_{jk}\df\frac1{(\mu\times\nu)(\mS_j\times \mT_k)}\int_{S_j\times T_k}\Phi\,d(\mu\times\nu).
%$$
%Тогда $\big\|\{c_{jk}\}\big\|_\frak M\le\|\Phi\|_{\frak M(\mu,\nu)}$.
%\end{lem}
%
%
%Рассмотрим ингегральный оператор с ядром $k(s,t)$ таким, что $k(s,t)=u_{jk}$ при $(s,t)\in\mS_j\times\mT_k$.
%
%Тогда 
%$$
%\int_{\mT}k(s,t)x(t)\,d\mu(t)=\sum_{k=1}^m u_{jk}\int_{\mT_k}x\,d\nu,\quad \text{если}\quad s\in\mS_j.
%$$
%

\begin{lem}  
\label{PQF0}
Let $\Phi\in\frak M(\mu,\nu)$. Then $(\mathcal P\odot\mathcal Q)\Phi\in\frak M(\mu,\nu)$
and $\|(\mathcal P\odot\mathcal Q)\Phi\|_{\frak M(\mu,\nu)}\le\|\Phi\|_{\frak M(\mu,\nu)}$.
\end{lem}

\Pf Using Lemmata \ref{msup0} and \ref{PQK0}, we obtain
\begin{multline*}
\|(\mathcal P\odot\mathcal Q)\Phi\|_{\frak M(\mu,\nu)}=
\sup\{\|\mathcal I_{((P\odot\mathcal Q)\Phi)h}\|_{\mathcal B(L^2(\nu),L^2(\mu))}:
\|\mathcal I_h\|_{\mathcal B(L^2(\nu),L^2(\mu))}\le1\}\\[.2cm]
=\sup\{\|\mathcal I_{(P\odot\mathcal Q)(\Phi h)}\|_{\mathcal B(L^2(\nu),L^2(\mu))}:(P\odot\mathcal Q)h=h, 
\|\mathcal I_h\|_{\mathcal B(L^2(\nu),L^2(\mu))}\le1\}\\[.2cm]
\le\sup\{\|\mathcal I_{\Phi h}\|_{\mathcal B(L^2(\nu),L^2(\mu))}: 
\|\mathcal I_h\|_{\mathcal B(L^2(\nu),L^2(\mu))}\le1\}=\| \Phi\|_{\frak M(\mu,\nu)}. \quad\bl
\end{multline*}

\begin{lem}  
\label{mlem0}
Let $\Phi\in\frak M(\mu,\nu)$. Then for an arbitrary positive $\e$, there exist a measurable function 
$\Psi$ in $\frak M(S,T)$
such that $\|\Phi-\Psi\|_{L^2(\mu\times\nu)}<\e$ and
$\|\Psi\|_{\frak M(S\times T)}=\|\Psi\|_{\frak M(\mu,\nu)}\le\|\Phi\|_{\frak M(S,T)}$.
\end{lem}

\Pf The linear span of the family $\{\chi_{E\times F}\}_{E\in\frak A, F\in\frak B}$ of characteristic functions is dense $L^2(\mu\times\nu)$. Consequently, we can select finite measurable  partitions $\{S_j\}_{j=1}^m$ and  
$\{T_k\}_{k=1}^n$ so that
 $$
 \inf\{\|\Phi-K\|_{L^2(\mu\times\nu)}: ~K\in L^2(\mu\times\nu),~(P\odot\mathcal Q)K=K\}<\e.
 $$
 Put $\Psi\df (P\odot\mathcal Q)\Phi$. Clearly, $\Psi\in\frak M(S,T)$ and $\|\Phi-\Psi\|_{L^2(\mu\times\nu)}<\e$.
 Everything else follows from from Lemmata \ref{mmatr0} and \ref{PQF0}. $\bl$

% \begin{thm} Пусть $\Phi\in\frak M(\mu,\nu)$. Тогда существует функция $\Phi_0\in\frak M(S\times T)$
%такая, что $\Phi_0=\Phi$ почти всюду по отношению к мере $\mu\times\nu$ и 
%$\|\Phi_0\|_{\frak M(S\times T)}=\|\Phi\|_{\frak M(\mu,\nu)}$.
%\end{thm}
%
%\Pf  Заметим сначала, что достаточно доказать только
%неравенство 
%\bay
%\label{Fi0}
%\|\Phi_0\|_{\frak M(S\times T)}\le\|\Phi\|_{\frak M(\mu,\nu)}.
%\ey
%
%Применяя лемму \ref{mlem0} для $\e=2^{-n}$, где $n\ge1$, получим последовательность измеримых
%функций $\{\Phi_n\}_{n=1}^\be$ из пространства $\frak M(S,T)$ такую, что
%$\|\Phi_n-\Phi\|_{L^2(\mu\times\nu)}<2^{-n}$ и
%$\|\Phi_n\|_{\frak M(S\times T)}\le\|\Phi\|_{\frak M(\mu,\nu)}$ при всех $n\ge1$.
%Из первого неравенства следует, что $\lim\limits_{n\to\be}\Phi_n=\Phi$ почти всюду.
%Остаётся в качестве $\Phi_0$ взять предельную точку последовательности  $\{\Phi_n\}_{n=1}^\be$ 
%в слабой топологии $\s(\ell^\be(S\times T),\ell^1(S\times T))$ пространства $\ell^\be(S\times T)$ . $\bl$

\

\section{\bf The proof of the main part of Theorem \ref{pmtn}}

\

In this section we assume that $\mu$ and $\nu$ are $\s$-finite measures; as we have already observed, the case of $\s$-finite measures can be reduced to the case of finite measures, and so in proofs we may assume the measures finite.

\begin{thm} 
\label{Phizero}
Let $\Phi\in\frak M(\mu,\nu)$. Then there exits a function $\Phi_0$ in $\frak M(S\times T)$
such that $\Phi_0=\Phi$ almost everywhere with respect to $\mu\times\nu$ and
$\|\Phi_0\|_{\frak M(S\times T)}=\|\Phi\|_{\frak M(\mu,\nu)}$.
\end{thm}

\Pf  Applying Lemma \ref{mlem0} for $\e=2^{-n}$, where $n\ge1$, we obtain a sequence of measurable
functions $\{\Phi_n\}_{n=1}^\be$ in $\frak M(S,T)$ such that
$\|\Phi_n-\Phi\|_{L^2(\mu\times\nu)}<2^{-n}$ and
$\|\Phi_n\|_{\frak M(S\times T)}\le\|\Phi\|_{\frak M(\mu,\nu)}$ for all $n\ge1$.
It follow from the first inequality that $\lim\limits_{n\to\be}\Phi_n=\Phi$ almost everywhere.
We take for $\Phi_0$ a limit point of the sequence $\{\Phi_n\}_{n=1}^\be$ 
in the weak topology $\s(\ell^\be(S\times T),\ell^1(S\times T))$ of the space $\ell^\be(S\times T)$. 
Now it i clear that $\|\Phi_0\|_{\frak M(S\times T)}\le\|\Phi\|_{\frak M(\mu,\nu)}$. It remains to observe that
$\|\Phi\|_{\frak M(\mu,\nu)}=\|\Phi_0\|_{\frak M(\mu,\nu)}$ because $\Phi_0=\Phi$ almost everywhere with respect to the measure $\mu\times\nu$ and
$\|\Phi_0\|_{\frak M(\mu,\nu)}\le\|\Phi_0\|_{\frak M(S\times T)}$ by Corollary \ref{STmunu0}.
$\bl$

\medskip

Let us complete now the proof of Theorem \ref{pmtn}.

\medskip

{\bf The proof of the implication (i)$\Longrightarrow$(iii) of Theorem \ref{pmtn}.} It follows from Theorem \ref{Phizero} that without loss of generality we may assume that $\Phi\in\frak M(S\times T)$
and $\|\Phi\|_{\frak M(S\times T)}=\|\Phi\|_{\frak M(\mu\times\nu)}$. It remains to apply Lemma \ref{STmunu}. $\bl$

\

\section{\bf The main theorem for Schur multipliers with respect to spectral measures}
\label{osnova}

\

In this section we obtain the main theorem for Schur multipliers with respect to spectral measures
that is going to be deduced from the results already obtained in this paper.

\begin{thm}
\label{osnspmery}
Let $E_1$ and $E_2$ be spectral measures in a separable Hilbert space and
$\Phi\in\fM(E_1,E_2)$. Then
$$
\|\Phi\|_{\fM(E_1,E_2)}=\|\Phi\|_{L^\be(E_1)\,\otimes_{\rm h}\!L^\be(E_2)}.
$$
\end{thm}

\Pf The result follows  from Theorem \ref{pmtn} and  from the criterion of the membership 
in the space of Schur multipliers given in \S\:\ref{Dvoi}. $\bl$

\medskip

{\bf Remark.} It is seen from Theorem  \ref{pmtn} that in the definition of the norm in the Haagerup tensor product given in \S\;\ref{Dvoi} the infimum in \rf{normaHaag} is attained.

\
 
\begin{footnotesize}
 
\noindent
\begin{tabular}{p{8cm}p{15cm}}
A.B. Aleksandrov & V.V. Peller \\
St.Petersburg State University & St.Petersburg State University \\
Universitetskaya nab., 7/9  & Universitetskaya nab., 7/9\\
199034 St.Petersburg, Russia & 199034 St.Petersburg, Russia \\
\\

St.Petersburg Department &St.Petersburg Department\\
Steklov Institute of Mathematics  &Steklov Institute of Mathematics  \\
Russian Academy of Sciences  & Russian Academy of Sciences \\
Fontanka 27, 191023 St.Petersburg &Fontanka 27, 191023 St.Petersburg\\
Russia&Russia\\
email: alex@pdmi.ras.ru& email: peller@math.msu.edu
\end{tabular}
\end{footnotesize}

\end{document}